\documentclass[aos,nameyear,dvips]{arximspdf}

\doi{10.1214/09-AOS754}
\volume{38}
\issue{4}
\pubyear{2010}
\firstpage{1937}
\lastpage{1952}

\makeatletter

\newtheorem{lemma}{Lemma}
\newtheorem{theorem}{Theorem}
\newtheorem{corollary}[theorem]{Corollary}

\newproclaim{remark}{Remark}

\makeatother

\begin{document}
\begin{frontmatter}

\title{Consistency of objective Bayes factors as the model dimension grows}
\runtitle{Consistency of objective Bayes factors}

\begin{aug}
\author[A]{\fnms{El{\'{\i}}as} \snm{Moreno}\thanksref{t1}\ead[label=e1]{emoreno@ugr.es}},
\author[B]{\fnms{F. Javier} \snm{Gir{\'{o}}n}\ead[label=e2]{giron@uma.es}} and
\author[C]{\fnms{George} \snm{Casella}\corref{}\thanksref{t2}\ead[label=e3]{casella@stat.ufl.edu}}
\runauthor{E. Moreno, F. J. Gir{\'{o}}n and G. Casella}
\affiliation{University of Granada, University of M{\'{a}}laga and University of Florida}
\address[A]{E. Moreno\\
Department of Statistics\\
University of Granada\\
18071, Granada\\
Spain\\
\printead{e1}} %adresu isvedimo komanda gale!
\address[B]{F. J. Gir{\'{o}}n\\
Department of Statistics\\
University of M\'{a}laga\\
Malaga\\
Spain\\
\printead{e2}}
\address[C]{G. Casella\\
Department of Statistics\\
University of Florida\\
Gainesville, Florida 32611\\
USA\\
\printead{e3}}
\end{aug}

\thankstext{t1}{Supported by Ministerio de Ciencia y Tecnolog{\'{\i}}a, Grant SEJ-65200
and Junta de Andaluc{\'{\i}}a Grant SEJ-02814.}

\thankstext{t2}{Supported by NSF Grants DMS-04-05543, DMS-0631632 and SES-0631588.}

\pdfauthor{Elias Moreno, F. Javier Giron, George Casella}

% HISTORY:
\received{\smonth{3} \syear{2009}}
\revised{\smonth{7} \syear{2009}}

% ABSTRACT
\begin{abstract}
In the class of normal regression models with a finite number of
regressors, and for a wide class of prior distributions, a Bayesian
model selection procedure based on the Bayes factor is consistent
[Casella and Moreno \textit{J.~Amer. Statist. Assoc.}
\textbf{104} (2009) 1261--1271]. However, in models where the number of parameters increases as
the sample size increases, properties of the Bayes factor are not
totally understood. Here we study consistency of the Bayes factors for
nested normal linear models when the number of regressors increases
with the sample size. We pay attention to two successful tools for
model selection [Schwarz \textit{Ann. Statist.} \textbf{6} (1978)
461--464] approximation to the Bayes factor, and the Bayes factor for
intrinsic priors [Berger and Pericchi \textit{J. Amer. Statist.
Assoc.} \textbf{91} (1996) 109--122, Moreno, Bertolino and Racugno
\textit{J. Amer. Statist. Assoc.} \textbf{93} (1998) 1451--1460].

We find that the the Schwarz approximation and the Bayes factor for
intrinsic priors are consistent when the rate of growth of the dimension of
the bigger model is $O(n^{b})$ for $b<1$. When $b=1$ the Schwarz
approximation is always inconsistent under the alternative while the Bayes
factor for intrinsic priors is consistent except for a small set of
alternative models which is characterized.
\end{abstract}

% KEYWORDS

\end{frontmatter}

%s1 ###
\section{Introduction}

Statistical methodology based on Bayes factors is particularly suitable for
dealing with multiple hypotheses testing problems when the dimension of the
parameter spaces varies across models. In such cases, Bayesian and
frequentist model selection procedures do not necessarily agree as the typically
\textit{ad hoc} dimension corrections of the different frequentist criteria
do not provide the same results as those automatically produced by the
Bayesian procedures which select models according to the parsimony
principle. For a recent discussion on the topic see Gir{\'{o}}n et al.
(\citeyear{Gironetal06}).

In the class of normal linear regression models, consistency of Bayesian
variable selection procedures, and, in particular, those using intrinsic
priors, has been recently established in Casella et al. (\citeyear{Casellaetal09}).
There it was shown that, under mild regularity conditions, when sampling
from a given submodel of a regression model with $p$ regressors, the
probability of selecting the true model tends to one as the sample size $n$
tends to infinity, and the probability of selecting any other submodel tends
to zero. It was also shown that the Schwarz (\citeyear{Schwarz78}) approximation is, in
spite of its simplicity, an accurate tool for selecting linear models when
there is a small number of parameters and the sample size is moderate or
large. Those results were obtained for a fixed number of regressors $p$, and
hence a finite number of models. Other forms of consistency of Bayes factors
for variable selection using Zellner's $g$-prior with several hyperpriors on $g$
have been recently studied by Liang et al. (\citeyear{Liangetal08}). These forms of
consistency include the consistency when $R^{2}$ tends to one, when $n$
tends to infinity, and consistency under prediction for squared error loss.

However, in some applications, the number of models increases with the sample
size. For instance, clustering is an interesting model selection problem
where the number of models increases as the sample size increases, and the
question is whether consistency of the Bayesian model selection procedure
based on intrinsic priors also holds in this latter context. Certainly, it
will not be possible to consistently estimate the parameters of the
underlying models, but we wonder whether consistently selecting the true
model is still possible.

When the number of parameters increases with the sample size, an
analysis of the consistency of several frequentist and Bayesian
approximation criteria for model selection in linear models, including
the Schwarz approximation, was given in Shao (\citeyear{Shao97}).
However, the results obtained by Shao (\citeyear{Shao97}) do not
coincide with ours, as the consistency notion used by Shao is not the
same as the one we use here. Shao defines a true model to be the
submodel minimizing the average squared prediction error, and
consistency of a model selection procedure means that the selected
model converges in probability to this model. We consider the true
model to be the one from which the observations are drawn. Of course,
the true model may not be in the class of models we are considering. In
this case consistency does not hold although many Bayesian model
selection procedures choose models in the class that are located as
close as possible to the true one where closeness is related to a
specific ``natural'' metric [Casella et al. (\citeyear{Casellaetal09})].

We examine consistency in linear models of both the Bayes factors for
intrinsic priors and the Schwarz approximation (BIC), when the dimension of
the parameter space of the models increases with the sample size. We find
that both the Bayes factor for the intrinsic priors and BIC are consistent
under the null; however, they might be inconsistent under some alternative
sampling models. The consistency depends on the rate of divergence of the
dimension of the null and the alternative linear models. Roughly speaking,
the BIC and the Bayes factor for intrinsic priors are consistent when the
rate of growth of the dimension of the full model goes to infinity as $%
O(n^{b})$ for any $b<1$.

When $b=1$, the BIC is always inconsistent under the alternative while
for the Bayes factor for intrinsic priors there is an inconsistency
region which is located in a small neighborhood of the null model. This
neighborhood is characterized in terms of a ``distance'' to the null
sampling model. In particular, for the case of the oneway ANOVA, where
$b=1$, the Bayes factor for intrinsic priors is not consistent for all
alternative models. This finding is apparently in contradiction with
the results of Berger, Ghosh and Mukhopadhyay (\citeyear{BGM03}) who
find ``\textit{that suitable Bayes factors will be consistent},'' and
hence induces an apparent paradox. However, in Section \ref{sec:berger}
we are able to resolve the apparent contradiction and find that the
consistency result in Berger, Ghosh and Mukhopadhyay (\citeyear{BGM03})
is obtained by using a normal prior centered at the null with variance
tending to zero, a situation that typically is not obtained by the
intrinsic priors. We also observe that consistency is obtained for this
problem for priors that degenerate to a point mass.

The rest of the paper is organized as follows. In Section
\ref{sec:consist} we characterize the consistency of the BIC and the
Bayes factor for intrinsic priors for the usual linear regression model
for $b<1$, demonstrate the inconsistency of BIC for $b=1$ and
characterize the small inconsistency region for the Bayes for intrinsic
priors for $b=1$. Section \ref{sec:app} presents some models where the
results of Section \ref{sec:consist} apply, and Section
\ref{sec:berger} resolves the apparent paradox with the results of
Berger, Ghosh and Mukhopadhyay (\citeyear{BGM03}). Section
\ref{sec:disc} provides a short, concluding discussion, and there is an
Appendix with some technical material.

%s2 ###
\section{Consistency in linear models}
\label{sec:consist}

In this section we give, for normal linear regression models with parameters
increasing with the sample size, conditions under which the Bayes factor for
intrinsic priors asymptotically selects the correct model. The finding is
that the Bayes factor may not be a consistent model selector for all
parameter values, depending on the rate of divergence of the dimension of
models. When there is an inconsistency region, it is characterized in terms
of a ``distance'' from the alternative to
the null model. We also show that the BIC model selector is inconsistent
when sampling from the full model if $b=1$.

Let $\mathbf{y}=(y_{1},\ldots,y_{n})^{{\prime }}$ be a vector of
independent responses, $\mathbf{X}_{p}$ a design matrix of dimension $%
n\times p$, where $p$ is the number of explanatory variables, and let $%
\mathbf{X}_{i}$ denote a submatrix of $\mathbf{X}_{p}$ whose dimensions are $%
n\times i$. We compare the reduced sampling model $N(\mathbf{y}|\mathbf{X}%
_{i}\mathbf{\alpha }_{i},\sigma _{i}^{2} \mathbf{I}_{n})$, and the full
model $N(\mathbf{y}|\mathbf{X}_{p}\mathbf{\beta }_{p},\sigma _{p}^{2} %
\mathbf{I}_{n})$, where the regression parameter vectors $\mathbf{\alpha }%
_{i}=(\alpha _{1},\ldots,\alpha _{i})^{\prime }$, $\mathbf{\beta }_{p}=(\beta
_{1},\ldots,\beta _{p})^{\prime }$ and the variance errors $\sigma _{i}^{2}$, $%
\sigma _{p}^{2}$, are unknown. Note that the reduced model is nested in
the full model. The comparison is based on the Bayes factor of model
$M_{p}$ versus model $M_{i}$, and we remark that it cannot be computed
by using the reference prior, the usual objective priors, since they
are improper and hence defined up to an arbitrary multiplicative
constant. The so-called intrinsic priors that are given below solve
this difficulty [Berger and Pericchi (\citeyear{BP96}), Moreno,
Bertolino and Racugno (\citeyear{MBR98})]. These objective priors have
proven to behave very well for multiple testing problems [Casella and
Moreno (\citeyear{CM06})].

To derive the Bayes factor for intrinsic priors we start with the improper
reference priors $\pi ^{N}(\mathbf{\alpha }_{i},\sigma _{i})=c_{i}/\sigma
_{i}$ and $\pi ^{N}(\mathbf{\beta }_{p},\sigma _{p})=c_{p}/\sigma _{p}$
where $c_{i}$ and $c_{p}$ are arbitrary positive constants, so we consider
the following Bayesian models:%
\[
M_{i}\dvtx\biggl\{N(\mathbf{y}|\mathbf{X}_{i}\mathbf{\alpha }_{i},\sigma _{i}^{2} %
\mathbf{I}_{n}), \pi ^{N}(\mathbf{\alpha }_{i},\sigma _{i})=\frac{c_{i}}{%
\sigma _{i}}\biggr\}
\]%
and%
\[
M_{p}\dvtx\biggl\{N(\mathbf{y}|\mathbf{X}_{p}\mathbf{\beta }_{p},\sigma _{p}^{2} %
\mathbf{I}_{n}), \pi ^{N}(\mathbf{\beta }_{p},\sigma _{p})=\frac{c_{p}}{%
\sigma _{p}}\biggr\}.
\]%
Standard calculations [Moreno, Bertolino and Racugno
(\citeyear{MBR98}),
and Gir{\'{o}}n et al. (\citeyear{Gironetal06})] yield the following intrinsic prior for $(\mathbf{\beta }%
_{p},\sigma _{p})$, conditional on $(\mathbf{\alpha }_{i},\sigma _{i})$:
\[
\pi ^{I}(\mathbf{\beta }_{p},\sigma _{p}|\mathbf{\alpha }_{i},\sigma _{i})=%
\frac{2\sigma _{i}}{\pi (\sigma _{i}^{2}+\sigma _{p}^{2})}N_{p}\Bigl(\mathbf{%
\beta }_{p}|\mathbf{\tilde{\alpha}}_{i},(\sigma _{i}^{2}+\sigma _{p}^{2})%
\mathbf{W}_{p}^{-1}\Bigr),
\]%
where $\mathbf{\tilde{\alpha}}_{i}^{\prime }=(\mathbf{\alpha }_{i}^{\prime },%
\mathbf{0}^{\prime })$, $\mathbf{0}^{\prime }$ being the null vector of $p-i$
components, and $\mathbf{W}_{p}^{-1}=\frac{n}{p+1}(\mathbf{X}_{p}^{\prime }%
\mathbf{X}_{p})^{-1}$. Then,\vspace*{1pt} using the priors $\{\pi ^{N}(\mathbf{\alpha }%
_{i},\sigma _{i}),\pi ^{I}(\mathbf{\beta }_{p},\sigma _{p}|\mathbf{\alpha }%
_{i},\sigma _{i})\}$, the Bayes factor for comparing the model $M_{p}$ and $%
M_{i}$ is
%e1 ###
\begin{equation}  \label{eq:intrinsicBF}\quad
B_{pi}(\mathbf{y})=\frac{2}{\pi }(p+1)^{(p-i)/2}\int_{0}^{\pi /2}\frac{\sin
^{p-i}\varphi (n+(p+1)\sin ^{2}\varphi )^{(n-p)/2}}{(n\mathcal{B}%
_{ip}+(p+1)\sin ^{2}\varphi )^{(n-i)/2}}\,d\varphi,   %\label{eq:IntegralExpre}
\end{equation}%
where
\[
\mathcal{B}_{ip}=\frac{\mathit{RSS}_{p}}{\mathit{RSS}_{i}}=\frac{\mathbf{y}^{\prime }(\mathbf{I%
}_{n}-\mathbf{H}_{p})\mathbf{y}}{\mathbf{y}^{\prime }(\mathbf{I}_{n}-\mathbf{%
H}_{i})\mathbf{y}}
\]%
and $\mathbf{H}_{j}=\mathbf{X}_{j}(\mathbf{X}_{j}^{\prime }\mathbf{X}%
_{j})^{-1}\mathbf{X}_{j}^{\prime },j=i,p$, is the hat matrix.

We first extend the definition of distance from $M_{p}$ to $M_{i}$
regression models given in Casella et al. (\citeyear{Casellaetal09}) to account for
models for which the number of parameters and the sample size increase to
infinity and define the ``distance'' from $%
M_{p}$ to $M_{i}$ for a given sample size $n$ as%
\[
\delta _{pi}=\frac{1}{\sigma _{p}^{2}}\mathbf{\beta }_{p}^{\prime }\frac{%
\mathbf{X}_{p}^{\prime }(\mathbf{I}_{n}-\mathbf{H}_{i})\mathbf{X}_{p}}{n}%
\mathbf{\beta }_{p}.
\]

The asymptotic performance of the Schwarz approximation is given in
Theorem~\ref{th:BIC}, and that of the Bayes factor for the intrinsic priors is given
in Theorems \ref{th:consistency} and \ref{th:inconsistency}, and Corollary %
\ref{cor:consistency}. The proofs of these results depend on Lemma \ref{lem:sampling}.

In what follows $\lim_{n\rightarrow \infty }[M]$ $Z_{n}$ will denote
the limit in probability of the random sequence $\{Z_{n};n\geq 1\}$
under the assumption that we are sampling from model~$M$. This model
$M$ will have a fixed parameter sequence. Further, we will need to use
the doubly noncentral beta distribution with parameters $\upsilon
_{1}/2,\upsilon _{2}/2$ and noncentrality parameters $\lambda
_{1},\lambda _{2}$. One way to define this distribution is as follows.
If $Y_{1},Y_{2}$ are independent random variables with noncentral chi
square distributions $\chi ^{2}(y_{1}|\upsilon _{1},\lambda _{1})$ and
$\chi ^{2}(y_{2}|\upsilon _{2},\lambda _{2})$, respectively, then the
variable $X=Y_{1}{\Large /}(Y_{1}+Y_{2})$ follows the doubly noncentral
beta distribution $\operatorname{Be}(\upsilon _{1}/2,\upsilon _{2}/2;\lambda
_{1},\lambda _{2})$ [Johnson, Kotz and Balakrishnan (\citeyear{JKB95}),
page 502].
\begin{lemma}
\label{lem:sampling}

\begin{enumerate}
\item When sampling from model $M_{i}$ the distribution of the statistics $%
\mathcal{B}_{ip}$ is the beta distribution $\operatorname{Be}((n-p)/2,(p-i)/2)$, and when
sampling from model $M_{p}$ it is the noncentral beta distribution $%
\operatorname{Be}((n-p)/2,(p-i)/2);0,n\delta _{pi})$. %
%n\delta _{pi}=\frac{1}{2\sigma _{p}^{2}}\beta _{p}^{\prime }\mathbf{X}%
%_{p}^{\prime }(\mathbf{I}_{n}-\mathbf{H}_{i})\mathbf{X}_{p}\beta _{p}.

\item Let $\{X_{n},n\geq 1\}$ be a sequence of random variables
such that%
\[
X_{n}\sim \operatorname{Be}\biggl(\frac{n-p}{2},\frac{p-i}{2};0,n\delta _{pi}\biggr),\qquad n\geq 1.
\]%
If $i$ and $p$ vary with $n$ as $i=O(n^{a})$
and $p=O(n^{b})$, where $0\leq a\leq b\leq 1$, then:

\begin{longlist}[(iii)]
\item[(i)] If $a<b=1$, when sampling from model $M_{p}$
\[
\lim_{n\rightarrow \infty } [M_{p}] X_{n}=\frac{1-1/r}{\delta +1},
\]%
where the constant $r$ satisfies $r=\lim_{p\rightarrow \infty
}n/p>1 $, and $\delta =\lim_{n\rightarrow \infty }\delta _{pi}$.

\item[(ii)] If $a=b=1$, so there exist two positive constants
such that $r=\lim_{p\rightarrow \infty }n/p>1$ and
$s=\lim_{p\rightarrow \infty }n/i>1$, we have
\[
\lim_{n\rightarrow \infty } [M_{p}] X_{n}=\frac{1-1/r}{1+\delta -1/s}.
\]

\item[(iii)] If $b<1$,
\[
\lim_{n\rightarrow \infty } [M_{p}] X_{n}=\frac{1}{1+\delta }.
\]
\end{longlist}
\end{enumerate}
\end{lemma}
\begin{pf}
See the \hyperref[app:consistencyproofs]{Appendix}.
\end{pf}

%s2.1 ###
\subsection{Inconsistency of BIC}

In this linear model setting we now prove that the Schwarz approximation for
comparing $M_{p}$ against $M_{i}$ is inconsistent when sampling from $M_{p}$
under certain conditions as first noticed by Stone (\citeyear{Stone79}) in a special case.
\begin{theorem}
\label{th:BIC} For comparing model $M_{p}$ to model $M_{i}$, where $M_{i}$
is nested in~$M_{p}$, and $i=O(n^{a})$ and $p=O(n^{b})$, if $0<a\leq b<1$,
the Schwarz approximation,
\[
S_{pi}(\mathbf{y})=\exp  \biggl\{ \frac{i-p}{2}\log n-\frac{n}{2}\log \mathcal{%
B}_{ip} \biggr\}
\]%
is consistent under the null and the alternative. However, if  $b=1$ it is
inconsistent under any alternative model $M_{p}$ provided that $%
\lim_{n\rightarrow \infty }\delta _{pi}>0$.
\end{theorem}
\begin{pf}
Consistency under the null for both cases follows from part 1 of Lemma \ref%
{lem:sampling}. For $b<1$, we notice that the leading term of the Bayes
factor is the one involving the statistic $\mathcal{B}_{ip}(\mathbf{y)}$,
but from part (iii) of Lemma \ref{lem:sampling} the limit of the sequence $%
\mathcal{B}_{ip}(\mathbf{y)}$ is a number strictly smaller than 1, and,
therefore,
\[
\lim_{n\rightarrow \infty } [M_{p}] S_{pi}(\mathbf{y})=\infty .
\]

On the other hand, if $b=1$, then $p=n/r$ and $i=n/s$, where $r$ is a
positive number greater than $1$ and $s$ is a number greater than $r$, the
leading term of the exponent of the Schwarz approximation is now the first
one which is strictly negative. Therefore,
\[
\lim_{n\rightarrow \infty } [M_{p}] S_{pi}(\mathbf{y})=0
\]%
and the proof is complete.
\end{pf}

%s2.2 ###
\subsection{Consistency of the Bayes factors for intrinsic priors}

We now characterize the consistency of the Bayes factor for intrinsic
priors, first assuming that both $p$ and $n$ increase at the same rate; that
is, $r=\lim_{n\rightarrow \infty,p\rightarrow \infty }n/p$, is a strictly
positive number. We further assume that the limit of the distance $\delta
_{pi}$ is finite when $p$ and $n$ tend to infinity, and $i$ is either finite
or increases to infinity at a lower rate than $n$. Note that in this theorem
the constant $b=1$.
\begin{theorem}
\label{th:consistency} Suppose that, as the sample size increases,
models increase their number of parameters with rate $i=O(n^{a})$
and $p=O(n)$, where $0\leq a<1$, and $r=\lim_{n
,p\rightarrow \infty }n/p> 1$.

\begin{enumerate}
\item When sampling from the simpler model $M_{i}$, $\lim_{n\rightarrow
\infty }[M_{i}] B_{pi}(\mathbf{y})=0$.

\item When sampling from the alternative model $M_{p}$ there exists a
function $\delta (r)$ such that
%e2 ###
\begin{equation}\label{eq:star}
\lim_{n\rightarrow \infty }[M_{p}] B_{pi}(\mathbf{y})= \cases{
\infty, & \quad  if  $\displaystyle\lim_{n\rightarrow \infty }\delta
_{pi}>\delta (r)$,\cr
$0$, &\quad if  $\displaystyle\lim_{n\rightarrow \infty }\delta _{pi}<\delta
(r)$.}
\end{equation}
\end{enumerate}

Further, this function has the simple expression %
%e3 ###
\begin{equation}  \label{eq:simple}
\delta (r)=\frac{r-1}{(r+1)^{(r-1)/r}-1}-1
\end{equation}%
and is a decreasing convex function such that $\lim_{r\rightarrow
\infty }\delta (r)=0$.
\end{theorem}
\begin{pf}
We first prove consistency of $B_{pi}(\mathbf{y})$ under the simpler model $%
M_{i}$. The Bayes factor $B_{pi}$ in (\ref{eq:intrinsicBF}) can be written
as
\[
B_{pi}(\mathbf{y})=\frac{2}{\pi }\int_{0}^{\pi /2} \biggl( 1+\frac{n}{%
(p+1)\sin ^{2}\varphi } \biggr) ^{(n-p)/2} \biggl( 1+\frac{n\mathcal{B}_{ip}}{%
(p+1)\sin ^{2}\varphi } \biggr) ^{-(n-i)/2}\, d\varphi .
\]%
From Lemma \ref{lem:sampling} it follows that
\[
\lim_{p\rightarrow \infty } [M_{i}] \mathcal{B}_{ip}=\frac{r-1}{r}
\]%
and, replacing $n$ by $p r$, the Bayes factor for large $p$ can be
approximated by
\[
B_{pi}(\mathbf{y})\approx \frac{2}{\pi }\int_{0}^{\pi /2} \biggl( 1+\frac{r}{%
\sin ^{2}\varphi } \biggr) ^{p(r-1)/2} \biggl( 1+\frac{r-1}{\sin ^{2}\varphi }%
 \biggr) ^{(i-pr)/2}\, d\varphi .
\]%
As the integrand is a monotonic increasing function of the angle $\varphi $,
the sup is attained at $\varphi =\pi /2$, and, therefore, an upper bound on
the integrand is $(1+r)^{p(r-1)/2} r^{(i-pr)/2}$. Then, for large $p$, an
upper bound for the Bayes factor is
\[
B_{pi}(\mathbf{y})< \biggl[ \frac{(1+r)^{r-1}}{r^{r}} \biggr] ^{p/2} r^{i/2}.
\]%
As the function of $r$ enclosed in square brackets is strictly smaller than $%
1$ for $r>1$, and the rate of growth of $i$ is strictly smaller than that of
$p$, it follows that
\[
\lim_{p\rightarrow \infty }  \biggl[ \frac{(1+r)^{r-1}}{r^{r}} \biggr]
^{p/2} r^{i/2}=0
\]%
for all $r>1$, thus proving consistency of the Bayes factor for the
intrinsic prior under the reduced model $M_{i}$.

Consistency under the full model $M_{p}$ is established as follows. From
Lemma \ref{lem:sampling}, the limiting distribution of the statistics $%
\mathcal{B}_{ip}$ under $M_{p}$ is
\[
\lim_{p\rightarrow \infty } [M_{p}] \mathcal{B}_{ip}=\frac{1-1/r}{\delta
+1},
\]%
where $\delta $ is the limit of the the ``distance'' from the full model to
the reduced one, which only depends on the limiting behavior of the
parameters of the full model;
%$\mathbf{\beta }_{p}$ when and $\sigma _{p}^{2}$,
that is,
\[
\delta =\lim_{p\rightarrow \infty }\delta _{pi}=\lim_{p\rightarrow \infty }%
\frac{1}{\sigma _{p}^{2}}\beta _{p}^{\prime }\frac{\mathbf{X}_{p}^{\prime }(%
\mathbf{I}_{n}-\mathbf{H}_{i})\mathbf{X}_{p}}{p r}\beta _{p}.
\]%
Therefore, the Bayes factor $B_{pi}(\mathbf{y})$ for large values of $p$ can
be approximated by
\[
B_{pi}(\mathbf{y})\approx \frac{2}{\pi }\int_{0}^{\pi /2} \biggl( 1+\frac{r}{%
\sin ^{2}\varphi } \biggr) ^{p(r-1)/2} \biggl( 1+\frac{r-1}{(1+\delta )\sin
^{2}\varphi } \biggr) ^{(i-pr)/2} \,d\varphi .
\]%
We look at two cases, depending on the values of the parameter $\delta $.

For $\delta >1$, the Bayes factor is an increasing convex function of $p$
and this implies that the Bayes factor is always consistent.

For $\delta \leq 1$, the argument proceeds as follows. As the integrand is a
continuous increasing function of $\varphi $ for all $r$, $\delta $ and $p$,
then by the mean value theorem, there exists a unique value of $\varphi _{0}$%
, say $0\leq \varphi _{0}(r,p,\delta )\leq \pi /2$, such that for large $p$
the Bayes factor is approximated by
\[
B_{pi}(\mathbf{y})\approx  \biggl( 1+\frac{r}{\sin ^{2}\varphi _{0}(r,p,\delta
)} \biggr) ^{p(r-1)/2} \biggl( 1+\frac{r-1}{(1+\delta )\sin ^{2}\varphi
_{0}(r,p,\delta )} \biggr) ^{(i-pr)/2}.
\]%
The limit of the sequence $\{\varphi _{0}(r,p,\delta ),p\geq 1\}$ is seen to
be equal to $\pi /2$ for all~$r$, and $\delta \leq 1$, Thus, for large
values of $p$, recalling that $i=o(p^{b})$, we can further approximate the
Bayes factor by%
%e4 ###
\begin{equation}  \label{eq:BFlimit}
B_{pi}(\mathbf{y})\approx  \biggl[ (1+r)^{r-1} \biggl( 1+\frac{r-1}{1+\delta }%
 \biggr) ^{-r}  \biggr] ^{p/2}.
\end{equation}
(It can be checked numerically that even for moderate values of $p$ this
approximation is very accurate.) Note that when the expression in square
brackets is greater than $1$ consistency holds, and when smaller than $1$
the Bayes factor is inconsistent. The root of the equation%
\[
(1+r)^{r-1} \biggl( 1+\frac{r-1}{1+\delta } \biggr) ^{-r}=1,
\]%
is $\delta (r)$ of (\ref{eq:simple}), %\begin{equation}
proving the theorem.
\end{pf}

We remark that the function $\delta (r)$ only depends on the $\lim $ $n/p=r$%
. In addition to the limiting value as $r \rightarrow \infty$, we also have $%
\lim_{r \rightarrow 0}\delta(r)=(e-1)^{-1}$ and $\lim_{r \rightarrow
1}\delta(r)=[log(2)]^{-1}-1$. Notice that the case of equality in the limit (%
\ref{eq:star}) is not covered by the theorem. It happens that, in this
case, we cannot make a specific conclusion as there will be parameter values
for which there is, and there is not, consistency.

Theorem \ref{th:consistency} covers the case in which the dimension of the
parameter space grows at a rate strictly smaller than that of the sample
space. However, it does not cover the case where the dimension of the null
and the alternative space grow at the same rate as the sample size. This
case is covered in Theorem \ref{th:inconsistency}.
\begin{theorem}
\label{th:inconsistency} Suppose that, as the sample size increases,
the rates the models increase their number of parameters are $i=O(n) $
and $p=O(n)$, and there exists positive constants $r $
and $s$ such that $r=\lim_{n,p\rightarrow \infty }n/p$
and $s=\lim_{n,i\rightarrow \infty }n/i \ge 1$.

\begin{enumerate}
\item When sampling from the simpler model $M_{i}$, $\lim_{n\rightarrow
\infty }[M_{i}] B_{pi}(\mathbf{y})=0$.

\item When sampling from the alternative model $M_{p}$, there exists
a function $\delta (r,s)$ such that
\[
\lim_{n\rightarrow \infty }[M_{p}] B_{pi}(\mathbf{y})= \cases{
\infty , &\quad  if  $\displaystyle\lim_{n\rightarrow \infty }\delta
_{pi}>\delta (r,s)$,\cr
$0$, &\quad if  $\displaystyle\lim_{n\rightarrow \infty }\delta _{pi}<\delta
(r,s)$.}
\]

This function has the following simple explicit form: %
%e5 ###
\begin{equation}  \label{eq:simple2}
\delta (r,s)=\frac{r-1}{(r+1)^{{s(r-1)}/({r(s-1)})}-1}-1+\frac{1}{s}
\end{equation}
and it is a bounded decreasing convex function in $r$ for
fixed $s$ with $\delta (r,s)\leq 1/\log 2-1$ for all
$s>r>1$, and $\lim_{r\rightarrow \infty }\delta (r,s)=0$
for all $s$. Further, $\lim_{s\rightarrow \infty }\delta
(r,s)=\delta (r)$ of (\ref{eq:simple}).
\end{enumerate}
\end{theorem}
\begin{pf}
To prove consistency under the simple model $M_{i}$, from Lemma \ref%
{lem:sampling} it follows that
\[
\lim_{p\rightarrow \infty } [M_{i}] \mathcal{B}_{ip}=\frac{s(r-1)}{r(s-1)}
\]%
and, replacing $n$ by $p r$, and $i/s=pr/s$, the Bayes factor for large $p $
can be approximated by
%e6 ###
\begin{eqnarray}
&&B_{pi}(\mathbf{y})\approx \frac{2}{\pi }\int_{0}^{\pi /2} \biggl( 1+\frac{r}{%
\sin ^{2}\varphi } \biggr) ^{p(r-1)/2}\nonumber\\[-8pt]\\[-8pt]
&&\hspace*{77.1pt}{}\times \biggl( 1+\frac{s(r-1)}{(s-1)\sin
^{2}\varphi } \biggr) ^{-{pr(s-1)}/({2s})} \,d\varphi.\nonumber
\end{eqnarray}

As the integrand is a monotonic increasing function of the angle $\varphi $,
the supremum is attained at $\varphi =\pi /2$, and thus an upper bound of the
integrand is
\[
(1+r)^{p(r-1)/2} \biggl( 1+\frac{s(r-1)}{(s-1)} \biggr) ^{-{pr(s-1)}/({2s})}.
\]

Then, for large $p$, the Bayes factor is bounded from above by
\[
B_{pi}(\mathbf{y})< \biggl[ (1+r)^{r-1} \biggl( \frac{rs-1}{s-1} \biggr) ^{{%
-r(s-1)}/{s}} \biggr] ^{p/2},
\]%
but as the function of $r$ and $s$ enclosed in square brackets is strictly
smaller than 1 for $s>r>1$, it follows that the limit of the upper bound of
the Bayes factor is 0 for all $s>r>1$ thus proving consistency of the Bayes
factor for the intrinsic prior under the reduced model $M_{i}$.

Consistency under the full model $M_{p}$ is proven in a similar way to that
of Theorem \ref{th:consistency}. From Lemma \ref{lem:sampling}, the limiting
distribution of the statistics $\mathcal{B}_{ip}$ under $M_{p}$ is now
\[
\lim_{p\rightarrow \infty } [M_{p}] \mathcal{B}_{ip}=\frac{1-1/r}{1+\delta
-1/s},
\]%
where $\delta $ is the same as in Theorem \ref{th:inconsistency}. Following the same course of
reasoning as in Theorem \ref{th:consistency}, we finally arrive at the following new
approximation for the Bayes factor for large values of $p$:
\[
B_{pi}(\mathbf{y})\approx  \biggl[ (1+r)^{r-1}  \biggl( 1+\frac{r-1}{1+\delta
-1/s} \biggr) ^{r(s-1)/s} \biggr] ^{p/2}.
\]

As the expression in square brackets does not depend on $p$, the limiting
behavior of the Bayes factor depends on whether this expression is less than
or greater than~$1$. Therefore, the new value of the boundary for
consistency-inconsistency, $\delta (r,s)$, is the root of the equation
\[
(1+r)^{r-1}  \biggl( 1+\frac{r-1}{1+\delta -1/s} \biggr) ^{r(s-1)/s}=1,
\]%
which is (\ref{eq:simple2}). %\begin{equation}
%the bound stated in the theorem.
This proves the theorem.
\end{pf}
\begin{remark}\label{remark1}
For all $s\geq r>1$, the function $\delta (r,s)$ is
bounded by a number smaller than $1$. Note also that if the rate of growth
of $M_{i}$ is smaller than that of~$M_{p}$, that is, $s\rightarrow \infty $,
then it is easy to show that $\lim_{s\rightarrow \infty }\delta (r,s)=\delta
(r)$.

An extension of Theorem \ref{th:inconsistency} to the case where models $%
M_{i}$ and $M_{p}$ grow at a slower rate than the sample size; that is, $%
i=O(n^{a})$ and $p=O(n^{b})$, where $0\leq a=b<1$, can be regarded as a
limiting case of the preceding theorem where both $r$ and $s$ go to
infinity. So, we have the following corollary.
\end{remark}
\begin{corollary}
\label{cor:consistency} For $i=O(n^{a})$ and $p=O(n^{b})$
and $0\leq a \leq b<1$, the Bayes factor for intrinsic
priors is consistent if  $\lim_{n\rightarrow \infty }\delta _{pi}>0$.
\end{corollary}

%This means, in particular, that when comparing models such that $i=O(n^{a})$
%and $p=O(n^{b})$ for\textit{ }$0\leq a\leq b<1$, the Bayes factor for
%intrinsic prior is consistent.

%s3 ###
\section{Applications}
\label{sec:app}

We look at some practical models for which the results of the preceding
section can be applied, including various ANOVA models, the multiple change
point problem, the clustering problem and spline regression.

In particular, the classical ANOVA problem will be illustrated in some
detail. For instance, we will see that for the one-way ANOVA, and by
extension any full factorial completely randomized design, the Bayes factor
for intrinsic priors is inconsistent in a region around the null. However,
reducing the ANOVA model by eliminating interaction terms recovers
consistency.
%We also discuss the ANOVA with known variances, mainly for comparing the behavior of the Bayes factor for
%intrinsic priors with that of the\textit{ }Bayes factor considered\textit{
%}by\textit{ }Berger \textit{et al.} (2003).
%Lastly, we show how our results apply to changepoint models, clustering, and spline regression.

%s3.1 ###
\subsection{Homoscedastic ANOVA}

There is a subtle difference between an ANOVA with a full model
specification (including all interactions) and one with a reduced model
specification, as it results in different asymptotic rates. We present the
results for balanced models with the same number of observations per cell,
but they can easily be extended to cover the unbalanced case.

%s3.1.1 ###
\subsubsection{Full model specification}

\label{sec:full}

We give a detailed development for the one-way ANOVA, and then show how the
results apply to full factorial designs. The null sampling model of the
homoscedastic one-way ANOVA, $M_{1}$,
%the observations $y_{ij}$ follow the normal distribution $N(y_{ij}|0,\tau
%^{2})$, and under the alternative sampling model $M_{p}$, the distribution
%is $N(y_{ij}|\mu _{i},\sigma ^{2})$, where $i=1,\ldots,p$, $j=1,\ldots,r$. Since
%the dimension of the null is $1$, and the dimension of the alternative is $%
%n+1$, as a consequence of theorem 3 there is an inconsistency region which
%is given by the set of alternative models such that \textit{ }$%
%More precise, we have that
%$\infty $\textit{,   if  }$\lim_{p\rightarrow \infty }\delta
%_{p1}>\delta (r)$, \\
%$0$\textit{,  if  }$\lim_{p\rightarrow \infty }\delta _{p1}<\delta
%(r)$.%
%The computation of the \textquotedblleft distance\textquotedblright $\delta
%_{p1}$ is, in general, a formidable task, but for this simple case the
%computation turns out to be easy. It can be shown that $\delta _{p1}$ is
%In an analogous way, the homoscedastic ANOVA
where it is assumed that the means are equal to an unknown $\mu $, can be
written as
\[
M_{1}\dvtx \biggl\{ N(\mathbf{y}|\mu \mathbf{1}_{n}, \tau ^{2} \mathbf{I}%
_{n}),\pi ^{N}(\mu,\tau )=\frac{c}{\tau } \biggr\}
\]%
and the alternative model as
\begin{eqnarray*}
&&M_{p}\dvtx \Biggl\{ N(\mathbf{y}|\mathbf{X}_{p}\mathbf{\mu }_{p},\sigma ^{2} %
\mathbf{I}_{n}),\pi ^{I}(\mathbf{\mu }_{p},\sigma |\mu,\tau )\\
&&\hspace*{7.4pt}\qquad=HC^{+}(\sigma
|\mu,\tau )\prod_{i=1}^{p}N\biggl(\mu _{i}|\mu,\frac{\tau ^{2}+\sigma ^{2}}{2}%
\biggr) \Biggr\},
\end{eqnarray*}
where $c$ is an arbitrary positive constant, $\mathbf{1}_{n}$ denotes a
vector of $n$ components containing $1$'s, $\mathbf{X}_{p}$ is an $n\times p$
matrix such that the first $r$ rows are equal to the unit vector $\mathbf{e}%
_{1}$, the next $r$ rows are equal to the unit vector $\mathbf{e}_{2}$ and
so on, so that the last $r$ rows are equal to the unit vector $\mathbf{e}%
_{p} $ where the unit vector $\mathbf{e}_{j}$ has coordinate $1$ at the $j$th
position, and $HC^{+}(\sigma |\mu,\tau )$ represents the half Cauchy
prior density of $\sigma $, conditional on $\mu ,\tau $, on the positive
part of the real line.

Since the dimension of $M_{1}$ is $2$ and the dimension of $M_{p}$ is $n+1$,
Theorem \ref{th:consistency} shows that there is an inconsistency region
given by those alternative models with $\lim_{p\rightarrow \infty }\delta
_{p1}<\delta (r)$ where $\delta (r)$ is given in (\ref{eq:simple}). Thus,
when sampling from $M_{p}$ we have that%
\[
\lim_{n\rightarrow \infty } [M_{p}] B_{p1}(\mathbf{y})= \cases{
\infty, & \quad if  $\displaystyle\lim_{p\rightarrow \infty }\delta _{p1}>\delta
(r)$, \cr
0, &\quad if  $\displaystyle\lim_{p\rightarrow \infty }\delta _{p1}<\delta
(r)$,}
\]%
where the distance $\delta _{p1}$ is given by
\[
\delta _{p1}=\frac{1}{n\sigma ^{2}}\mathbf{\mu }_{p}\biggl(\mathbf{I}_{n}-\frac{1}{%
n}\mathbf{1}_{n}\mathbf{1}_{n}\biggr)\mathbf{\mu }_{p}=\frac{1}{\sigma ^{2}}\frac{1%
}{p}\sum_{i=1}^{p}(\mu _{i}-\bar{\mu}_{p})^{2}.
\]%
%We note that  $\delta(r)$ decreases to $0$ quite rapidly (for example, $\delta(r) \approx .25$ for $r=5$), so the region of inconsistency is almost of no concern.

If we have a multiway completely randomized design the same results hold, as
such a design is equivalent to a one-way design. For example, suppose we have
a three-way full factorial with the model
%e7 ###
\begin{eqnarray}  \label{eq:threeway}
y_{ijk}&=&\mu_i+\tau_j+\gamma_k+(\mu\tau)_{ij}+(\mu\gamma)_{ik}+(\tau%
\gamma)_{jk}+(\mu\tau\gamma)_{ikj}+\varepsilon_{ijk},
\nonumber\\[-8pt]\\[-8pt]
\eqntext{i = 1, \ldots, I,   j = 1, \ldots, J,   k = 1, \ldots, K.}
\end{eqnarray}
The number of parameters (with identifiability restrictions) is $IJK$, and
thus we are again in the case of Theorem \ref{th:consistency} with $b=1$.
Any null hypothesis will result in a model $M_1$ with a reduced set of
parameters that will satisfy $a<b$ of Theorem \ref{th:consistency}. Thus
when sampling from the full model, the intrinsic Bayes procedure is
consistent only if $\delta _{p1}>\delta(r)$, where $\delta(r)$ is given in (%
\ref{eq:simple}), and, analogous to the one-way case, $\delta _{p1}$ is equal
to the sum of squares of the differences between the null model coefficients
and the full model coefficients. Extension to higher-order designs is
straightforward.

%s3.1.2 ###
\subsubsection{Reduced model specification}

In higher-order ANOVA models, it is often the case that some interaction
terms are not specified. In particular, if the highest order interaction is
not in the model, we can attain consistency of the intrinsic Bayes factor
over the entire parameter space. We illustrate this with the three-way model (%
\ref{eq:threeway}); the extension to higher-order models should be clear.

If we eliminate the term $(\mu \tau \gamma )_{ikj}$ from the model (\ref%
{eq:threeway}), then there are at most
\[
p=I+J+K+IJ+IK+JK
\]%
parameters in the full model $M_{2}$. Since there are $n=rIJKL$
observations, it immediately follows that
\[
p= \cases{
O(n), &\quad if $I$ or $J$ or $K\rightarrow \infty$, \cr
o(n), &\quad if $I$ and $J$ and $K\rightarrow \infty$.}
\]
So in the first case we can apply Theorem \ref{th:consistency}, and, similar
to the full model evaluation, there will be an inconsistency region.
However, in the second case, when all of $I$, $J$, and $K$ $\rightarrow
\infty $ we are in the case of Corollary \ref{cor:consistency}; there is no
inconsistency region and the Bayes factor for the intrinsic priors is
consistent in the entire parameter space.

%s3.2 ###
\subsection{Nested regression models}

Clustering, multiple change points and spline regression are examples of model
selection problems for which the dimension of the alternative models grows
at the same rate as the sample size $n$. Therefore, in the notation of the
preceding sections $b=1$, and hence the Schwarz approximation is
inconsistent, but the Bayes factor for intrinsic priors is consistent except
for a small region around the null model. Note that the null model in
clustering is the one cluster model, in the multiple change points problem
the null model is the no change model, and in spline regression the null
model is the model that specifies no knots.

%s4 ###
\section{Comparison with previous findings}
\label{sec:berger}

As we have seen in Section \ref{sec:full}, in the homoscedastic ANOVA there
is a region of inconsistency for the Bayes factor for intrinsic priors. This
result seems to be in contradiction with the finding by Berger, Ghosh and
Mukhopadhyay (\citeyear{BGM03}) who consider the Bayes factor for normal priors. The models they
compare are essentially%
%e8 ###
\begin{eqnarray}\label{eq:BergerCompare}
&&M_{1} \dvtx \Biggl\{
\prod _{i=1}^{p}\prod _{j=1}^{r}N(y_{ij}|0,1) \Biggr\} ,
\nonumber\\[-8pt]\\[-8pt]
&&M_{2} \dvtx \Biggl\{ \prod _{i=1}^{p}\prod _{j=1}^{r}N(y_{ij}|\mu
_{i},1), \pi ^{I}(\mathbf{\mu }_{p}|t)=\prod _{i=1}^{p}N(\mu
_{i}|0,2/t), t\geq 1 \Biggr\},   \nonumber
\end{eqnarray}%
where $\pi ^{I}(\mathbf{\mu }_{p}|t)$ is the intrinsic prior when a training
sample of size $t$ is considered. Observe that the hyperparameter $t$
controls the degree of concentration of the intrinsic priors around the
null, and it usually ranges from $1$ to $r$ so as to not exceed the
concentration of the likelihood of $\mu _{i}$ [for a discussion on the topic
see Casella and Moreno (\citeyear{CM09})].

For a given sample $\mathbf{y}=\{y_{ij}$, $j=1,\ldots,r$, $i=1,\ldots,p\}$, the
Bayes factor for intrinsic priors to compare the Bayesian model $M_{2}$
against $M_{1}$ is
\[
B_{21}(\mathbf{y}|t)= \biggl( \frac{t}{2r+t} \biggr) ^{p/2}\exp  \Biggl\{ \frac{%
r^{2}}{2r+t}\sum_{i=1}^{p}\bar{y}_{i}^{2} \Biggr\}
\]%
and it satisfies
%e9 ###
\begin{equation}\label{eq:MorenoCond}
\lim_{p\rightarrow \infty } [M_2] B_{21}(\mathbf{y}|t)= \cases{
\infty, &\quad if $\displaystyle\lim_{p\rightarrow \infty }\frac{1}{p}%
\sum_{i=1}^{p}\mu _{i}^{2}>R(t,r)$, \cr
0, &\quad if $\displaystyle\lim_{p\rightarrow \infty }\frac{1}{p}%
\sum_{i=1}^{p}\mu _{i}^{2}<R(t,r)$,}
\end{equation}%
where $R(t,r)=(2r+t)(2r^{2})^{-1}\ln [(2r+t)t^{-1}]-r^{-1}$, $1\leq t\leq r$.

As a curiosity we mention that the function $R(t,r)$ is related to the
function $\delta (r)$ of Theorem \ref{th:consistency} in the following way:
\[
R(2,r)<\delta (r)<R(1,r),\qquad r\geq 1.
\]

The Bayes factor for intrinsic priors is not consistent for all
possible alternative sampling models, and thus we cannot call it a
consistent model selector. For each $t$ and $r $, the inconsistency region
in the alternative parametric space will be denoted as
%e10 ###
\begin{equation}
C(t,r)=\Biggl\{\mathbf{\mu }\dvtx0<\lim_{p\rightarrow \infty }\frac{1}{p}%
\sum_{i=1}^{p}\mu _{i}^{2}<R(t,r)\Biggr\}.
\end{equation}%
We note that the bound $R(t,r)$ is a decreasing function in both arguments $%
t $ and $r$, and $\lim_{r\rightarrow \infty }R(t,r)=\lim_{t\rightarrow
\infty }R(t,r)=0$.

It turns out that for some extreme priors the Bayes factor is a consistent
model selector. We present two extreme cases: the first one where the prior
degenerates to a point mass, and the second one for intrinsic priors with
variances that tend to zero.

1. \textit{Simple null versus simple alternative}. As a
modification of (\ref{eq:BergerCompare}), suppose we want to choose between
\[
M_{1}\dvtx\prod _{i=1}^{p}\prod _{j=1}^{r}N(y_{ij}|0,1)%
\quad\mbox{and}\quad M_{2}\dvtx\prod _{i=1}^{p}\prod _{j=1}^{r}N(y_{ij}|\mu
_{i0},1),
\]%
where $\{\mu _{i0},i\geq 1\}$ is an arbitrary but specified sequence such
that\break $\lim_{p\rightarrow \infty }\sum \mu _{i0}^{2}/p>0$. Then, the Bayes
factor $B_{21}(\mathbf{y}_{p})$ satisfies\break $\lim_{p\rightarrow \infty } [M_{2}] B_{21}(\mathbf{y})=\infty $; that is,
the Bayes factor is consistent under the alternative. This simple result
means that when the prior distribution on the alternative degenerates to a
point mass, consistency of the corresponding Bayes factor holds.
%In other words, no uncertainty on the alternative parameter space yields
%consistency of the Bayes factor.

2. \textit{Mixture priors}. The presence of uncertainty in the alternative
models provokes the appearance of an inconsistency region $C(t,r)$. However,
in Berger, Ghosh and Mukhopadhyay (\citeyear{BGM03}), they use a continuous version of the
intrinsic prior above and augment $M_{2}$ by mixing the variance $1/t$ of
the $N(\mu _{i}|0,1/t)$ with a hyperprior density $g(t)$. (Special cases
they consider are to take $g$ to be either gamma or beta, yielding priors
that they refer to as Cauchy and Smooth Cauchy.) For these general mixture
priors they prove the following theorem.
\begin{theorem*}[{[Berger, Ghosh and Mukhopadhyay (\citeyear{BGM03}), Theorem 3.1]}]
For any
prior of the form
%e11 ###
\begin{equation}  \label{eq:BergerPrior}
\pi _{g}(\mathbf{\mu })=\int_{0}^{\infty }\frac{t^{p/2}}{(2\pi )^{p/2}}%
e^{-(t/2)\sum_{i}\mu _{i}^{2}}g(t)\,dt
\end{equation}%
with $g(t)$ having support on $(0,\infty )$, the
Bayes factor is consistent under $M_{1}$. Consistency under $M_{2}$
holds if
%e12 ###
\begin{equation}\label{eq:BergerCond}
\tau ^{2}=\lim_{p\rightarrow \infty }\frac{1}{p}\sum_{i}\mu _{i}^{2}>0.
\end{equation}
\end{theorem*}

How do we reconcile (\ref{eq:MorenoCond}) and (\ref{eq:BergerCond}), an
apparent paradox? To obtain consistency for any alternative sampling model,
we need the function in (\ref{eq:MorenoCond}) to be zero, but this only
occurs when $t$ goes to infinity because $r$ is fixed. Since the
inconsistency regions $\{C(t,r),t\geq 1\}$ form a monotone decreasing
sequence, the limit is $C_{\infty }(r)=\bigcap_{t=1}^{\infty }C(t,r)=\{\mu
_{i}\dvtx\lim_{p\rightarrow \infty }\frac{1}{p}\sum_{i}\mu _{i}^{2}=0\}$, a
point that does not belong to the alternative parameter space. In the above
theorem this is exactly what the prior $\pi _{g}(\mathbf{\mu })$ does by
incorporating priors with variance that tends to zero. (Something similar
produces the so-called Lindley paradox when testing that the mean of a
normal is zero; as the variance of the normal prior goes to zero less and
less prior mass is given to any neighborhood of the null.)

Certainly, if we mix values of $t$ from $1$ to $r<\infty $, for instance
mixing all the intrinsic priors, the intersection of these inconsistency
regions $C_{r}(r)=\bigcap_{t=1}^{r}C(t,r)$, is a nonempty set in the
alternative model space, and hence the inconsistency region does not
disappear. This is also noted by Berger, Ghosh and Mukhopadhyay in their Theorem
3.2. that we state here using our notation.
\begin{theorem*}[{[Berger, Ghosh and Mukhopadhyay (\citeyear{BGM03}), Theorem 3.2]}]
For any
prior of form (\ref{eq:BergerPrior}), with $g(t)$ being supported
on a finite interval $[0,1]$, and $r=1$, the Bayes Factor
is inconsistent under $M_{2}$ for $0<\tau ^{2}<2\log 2-1$.
\end{theorem*}

We note that here $t=2$ and $R(2,1)=2\log 2-1$.

%s5 ###
\section{Discussion}
\label{sec:disc}

In our previous work [Casella et al. (\citeyear{Casellaetal09})], where we looked at
consistency of Bayes factors for a fixed number of parameters, we found that
both the Bayes factor for intrinsic priors and the Schwarz approximation to
a Bayes factor had the same asymptotic behavior, and both were consistent.
In this paper we have derived the asymptotic behavior of the Bayes factor
for intrinsic priors and the Schwarz approximation when the dimension of the
model grows with the sample size, and we note an interesting dichotomy in
their performance. The Bayes factor for intrinsic priors and the Schwarz
approximation have very different asymptotic behavior for the usual case
where the dimension of the full model grows at the same rate as the sample
size with the Bayes factor for intrinsic priors clearly being the optimal
one.

\begin{table}
\tablewidth=200pt
\caption{}\label{table1}
\begin{tabular*}{\tablewidth}{@{\extracolsep{\fill}}ll@{}}
\hline
\textbf{Rate of divergence} & \textbf{Consistency region of} $\bolds{B_{pi}}\bolds{(}\mathbf{y}\bolds{)}$ \\
\hline
$0 < a=b=1$ & $M_{p}\dvtx\lim_{n\rightarrow \infty }\delta _{pi}>\delta (r,s)$
\\
$0\leq a<b=1$ & $M_{p}\dvtx\lim_{n\rightarrow \infty }\delta _{pi}>\delta (r)$
\\
$0\leq a\leq b<1$ & $M_{p}\dvtx\lim_{n\rightarrow \infty }\delta _{pi}>0$ \\
\hline
\end{tabular*}
\end{table}

We summarize the consistency regions of the Bayes factor for intrinsic
priors for different values of $a$ and $b$ in
Table \ref{table1}, and we extract the following recommendations. For models with $%
b<1$, the existence of very many parameters is not an inconvenience as far
as consistency is concerned. For models with $b=1$, there is a small
inconsistency region around the null defined by the function $\delta $ that
decreases rapidly as $r$ increases. It also follows that inconsistency is
the exception for the Bayes factor for intrinsic priors, while the rule is
consistency, and this gives credence to the Bayes factor for intrinsic
priors as a powerful objective tool for model selection.

\begin{appendix}\label{app:consistencyproofs}

\section*{\texorpdfstring{Appendix: Proof of Lemma
\protect\lowercase{\ref{lem:sampling}}}{Appendix: Proof of Lemma 1}}

Part 1 follows from Theorem 1 in Casella et al. (\citeyear{Casellaetal09}). To prove
part 2, we note that $X_{n}$ can be written as
\[
X_{n}= \biggl( 1+\frac{V_{n}}{W_{n}} \biggr) ^{-1},
\]%
where $V_{n}\sim (1/n)\chi _{p-i}^{2}(n\delta _{pi})$ and $W_{n}\sim
(1/n)\chi _{n-p}^{2}$. The means and variances of these random variables are
\[
\mathrm{E}(V_{n})=\delta _{pi}+\frac{p-i}{n}, \qquad \mathrm{E}(W_{n})=1-\frac{%
p}{n}
\]%
and
\[
\operatorname{Var}(V_{n})=\frac{4\delta
_{pi}}{n}+\frac{2(p-i)}{n^{2}},\qquad
\operatorname{Var}(W_{n})=\frac{2(n-p)}{n^{2}}.
\]%
From these expressions the three cases follow:

\begin{longlist}[(iii)]
\item[(i)] If $a<b=1$, then when sampling from model $M_{p}$
\[
V_n \rightarrow \delta+\frac{1}{r},  \qquad W_n \rightarrow
1-\frac{1}{r}\quad
\mbox{and}\quad X_n \rightarrow \frac{1-{1}/{r}}{1+\delta}.
\]

\item[(ii)] If $a=b=1$, then when sampling from model $M_{p}$
\[
V_n \rightarrow \delta+\frac{1}{r}-\frac{1}{s},  \qquad W_n \rightarrow 1-%
\frac{1}{r}\quad\mbox{and}\quad X_n \rightarrow \frac{1-{1}/{r}}{1+\delta-{%
1}/{s}}.
\]

\item[(iii)] If $b<1$, then when sampling from model $M_{p}$
\[
V_n \rightarrow \delta, \qquad W_n \rightarrow 1\quad\mbox{and}\quad X_n
\rightarrow \frac{1}{1+\delta}.
\]
\end{longlist}
\end{appendix}

\printaddresses

\end{document}